\documentclass[12pt]{article}

\usepackage{geometry} 
\usepackage{amsthm,amssymb}
\usepackage{amsmath,amsfonts,latexsym}
\usepackage{latexsym}
\usepackage{float}
\usepackage{xcolor}
\usepackage{url}
\usepackage{comment}

\def\titlerunning#1{\gdef\titrun{#1}}
\makeatletter
\def\author#1{\gdef\autrun{\def\and{\unskip, }#1}\gdef\@author{#1}}
\def\address#1{{\def\and{\\\hspace*{18pt}}\renewcommand{\thefootnote}{}%
\footnote {#1}}%
\markboth{\autrun}{\titrun}}
\makeatother
\def\email#1{\hspace*{4pt}{\em e-mail}: #1}

\newtheorem{thm}{Theorem}[section]
\newtheorem{prop}[thm]{Proposition}
\newtheorem{lemma}[thm]{Lemma}
\newtheorem{cor}[thm]{Corollary}
\theoremstyle{definition}
\newtheorem{rem}[thm]{Remark}

\newcommand{\We}{\mathrm{W}}

\newcommand{\F}{\mathbb{F}}

\begin{document}

\titlerunning{}
\title{LCD codes from weighing matrices}

\author{Dean Crnkovi\' c, Ronan Egan, B. G. Rodrigues and Andrea \v{S}vob }

\date{\vspace{-5ex}}

\maketitle

\address{D. Crnkovi\'{c}, A. \v{S}vob: Department of Mathematics, University of Rijeka, Croatia;
\email{\{deanc,asvob\}@math.uniri.hr}
\and
R. Egan: School of Mathematics, Statistics and Applied Mathematics, National University of Ireland, Galway;
\email{ronan.egan@nuigalway.ie}
\and
B. G. Rodrigues: School of Mathematics, Statistics and Computer Science, University of KwaZulu-Natal, Durban, South Africa;
\email{rodrigues@ukzn.ac.za}
}

\begin{abstract}
Linear codes with complementary duals are linear codes whose intersection with their duals are trivial, shortly named LCD codes. In this paper we outline a construction for LCD codes over finite fields of order $q$ using weighing matrices and their orbit matrices. The LCD codes constructed can be of any length dimension according to the choice of matrices used in their construction.  As a special case, LCD codes of length $2n$ and dimension $n$ are constructed which also have the property of being formally self-dual. Alternatively, under a condition depending on $q$ that the codes are not LCD, this method constructs self-dual codes.  To illustrate the method we construct LCD codes from weighing matrices, including the Paley conference matrices and Hadamard matrices. We also extend the construction to Hermitian LCD codes over the finite field of order 4. In addition, we propose a decoding algorithm that can be feasible for the LCD codes obtained from some of the given methods.
\end{abstract}

\bigskip



\section{Introduction}

In this paper we introduce a method for constructing linear codes with complementary duals, or LCD codes. LCD codes are linear codes whose intersection with their duals are trivial; they were introduced by Massey in \cite{massey} and have been widely applied in information protection, electronics and cryptography.  They provide an optimum linear coding solution for the two user binary adder channel.
LCD codes are asymptotically good; it is shown in \cite{sendrier} that they meet the asymptotic Gilbert-Varshamov bound.  Carlet and Guilley \cite{carlet-guilley} studied applications of binary LCD codes against side-channel attacks (SCA) and fault injection attacks (FIA). In their work, they presented several constructions of LCD codes and showed that non-binary LCD codes in characteristic 2 can be transformed into binary LCD codes by expansion.  For further recent research on the topic we refer the reader to \cite{ArayaHarada,doug, harada-so-lcd, sole}.  One of our main challenges is to construct LCD codes over finite fields that also have good error correcting properties. A result of Carlet et al. in \cite{carlet-et-al} gives that over any finite field of order $q > 3$, the existence of an
$[n,k,d]_{q}$ linear code implies the existence of an $[n,k,d]_{q}$ LCD code.  The limits to error correction may be more restrictive for binary and ternary LCD codes than for linear codes in general,
as LCD codes do satisfy extra constraints.

Formally self-dual codes, whose weight enumerators are invariant under MacWilliams relations, are also studied because they occasionally have stronger error correcting properties than self-dual codes may achieve.  They may also, in certain cases, contain within the codewords some combinatorial structures such as $t$-designs \cite{KennedyPless}.

Our primary aim in this paper is to introduce a direct construction of LCD codes using weighing matrices, $(r, \lambda)$-designs, and their orbit matrices.
There are many known constructions of weighing matrices, especially Hadamard matrices, and $(r, \lambda)$-designs. For example, for each prime power $q$, $q \equiv 3\ (mod\ 4)$,
there is a Paley type I Hadamard matrix of order $(q+1)$, and for each prime power $q$, $q \equiv 1\ (mod\ 4)$, there is a Paley type II Hadamard matrix of order  $2(q+1)$
and a weighing matrix of order $(q+1)$ called a conference matrix. Further, there are a lot of known Hadamard matrices and the famous Hadamard conjecture proposes that for every positive integer
$m$ there exists a Hadamard matrix of order $4m$ (see \cite{seberry}). In addition, every weighing matrix can lead us to various orbit matrices.
Therefore, the methods described in this paper can produce a large number of LCD codes.
We construct $[n,k]$ LCD codes for any length $n$ and dimension $k$ over most finite fields.  As a special case we construct $[2n,n]$ LCD codes which will also have the property of being formally self-dual. When the weighing matrix does not satisfy the given conditions for the constructed code to be an LCD code over $\F_{q}$, then the code will be self-dual.

When Massey \cite{massey} introduced the terminology for LCD codes, he also introduced a decoding method for LCD codes that involved a map $\varphi$ from the dual code
$C^\perp$ to $C$ where for $v \in C^\perp$, $\varphi(v)$ is the codeword in $C$ closest to $v$. We show how this can be done for some of the constructed codes, using a method that is feasible
for a small number of errors. The approach that we use to construct the map $\varphi(v)$ that is involved in a decoding algorithm is based on the approach used in \cite{LCD-Peisert}.

The paper is outlined as follows.  Section \ref{Preliminaries} provides the necessary definitions and notation used throughout.  In Section \ref{lcd_mat} the main construction of LCD codes using weighing matrices with various properties is outlined, then in Section \ref{lcd_orb_mat} these ideas are extended to orbit matrices of weighing matrices with respect to certain automorphism groups.  In the final section we extend the approach to the construction of Hermitian LCD codes over the finite field of 4 elements.  Examples of codes constructed are given throughout the paper.  These codes have been constructed and examined using Magma \cite{magma}.

\section{Preliminaries}\label{Preliminaries}

We adhere mostly to the standard coding theory notation of \cite{Huffman}, to which we refer for background reading. Throughout, we let $\F_{q}$ denote the finite field of order $q$ for some prime power
$q$.  We let $\dagger$ be the transposition acting on a field $\F$ such that $0^{\dagger} = 0$ and $x^{\dagger} = x^{-1}$ for all $x \in \F\setminus \{0\}$.
When $M = [m_{ij}]$ is an $n \times n$ matrix with entries in $\F$, we let $M^{*} = [m_{ji}^{\dagger}]$. This is to coincide with the standard use of $M^{*}$ to denote the conjugate transpose of a matrix over the complex field, when entries are either zero, or are of modulus $1$.  Clearly $M^{*} = M^{\top}$ when the entries are in the set $\{0,\pm 1\}$.
If $v = [v_{0},\ldots,v_{n-1}]$ is a vector, then $v^{*} = [v_{0}^{\dagger},\ldots,v_{n-1}^{\dagger}]$. Let $I_n$ and $J_{n}$ denote the $n \times n$ identity and all-ones matrix over some ring or field
which will be clear from the context.  If $W$ has entries in $\{0,\pm 1 \}$ and $WW^{*} = mI_{n}$ over $\mathbb{C}$, $W$ is a \emph{weighing matrix} $\We(n,m)$, and if $m=n$,
$W$ is a \emph{Hadamard matrix} $\mathrm{H}(n)$. In this paper, weighing matrices have entries in $\{0,\pm 1\}$, and are often embedded as submatrices into matrices with entries in $\F_{q}$. This is done with the understanding that the entries are interpreted as $0$, $1$ or $-1$ in $\F_{q}$ accordingly.

Let $\mathcal{V}$ be a set of $v$ points, and let $\mathcal{B}$ be a set of $b$ subsets of $\mathcal{V}$ known as blocks.  Suppose that for a fixed set of integers $K$, every block of $\mathcal{B}$ is a
$k$-subset of $\mathcal{V}$ with $k \in K$.  Then $(\mathcal{V},\mathcal{B})$ is a \emph{pairwise balanced design} $\mathrm{PBD}(v,K,\lambda)$ if every pair of distinct elements of $\mathcal{V}$ are both incident with precisely $\lambda$ blocks of $\mathcal{B}$.  If $K = \{k\}$ is a singleton set, then this is a $2$-\emph{design} or \emph{balanced incomplete block design} $\mathrm{BIBD}(v,k,\lambda)$.  Where $K$ is unspecified, the pair $(\mathcal{V},\mathcal{B})$ is an $(r, \lambda)$-design if each point of $\mathcal{V}$ is incident with precisely $r$ blocks of $\mathcal{B}$, and every distinct pair of elements are both incident with precisely $\lambda$ blocks of $\mathcal{B}$.  If $(\mathcal{V},\mathcal{B})$ is a BIBD, then each point is incident with a constant number of blocks, denoted by $r$, and called a replication number, and so is an $(r, \lambda)$-design.  If $B$ is the $v \times b$ point-by-block incidence matrix of a $(r,\lambda)$-design, then
$BB^{\top} = (r-\lambda)I_{b} + \lambda J_{b}$.  For more information on $(r,\lambda)$-designs we refer the reader to \cite{crc-rl-des}.

A \emph{$q$-ary linear code} $C$ of length $n$ and dimension $k$ is a $k$-dimensional subspace of $\F_{q}^{n}$.  We say the $k \times n$ matrix comprised of rows that span $C$ \emph{generates}, or is a \emph{generator matrix} of $C$. We refer only to linear codes in this paper.

The \emph{dual} code $C^\perp$ is the orthogonal complement under the standard inner product
$\langle\cdot \, ,\cdot\rangle$, i.e.\ $C^{\perp} = \{ v \in \F_{q}^n| \langle v,c\rangle=0 {\rm \ for\ all\ } c \in C \}$.  A linear code $C$ over $\F_q$ is called a Euclidean or classical LCD code if
$C\cap C^{\perp}=\{0\}$.  Usually we just write LCD code in this instance. It follows that the dual of an LCD code is an LCD code.
A code $C$ is \emph{self-orthogonal} if $C \subseteq C^\perp$ and \emph{self-dual} if equality is attained.

Given codewords $x=(x_1,...,x_n)$ and $y=(y_1,...,y_n)\in \mathbb{F}_q^n$,
the \emph{Hamming distance} between $x$ and $y$ is the number \ $d(x,y)=\left| \{ i  :  x_i \neq y_i \} \right| $.
The \emph{minimum distance} of the code $C$ is defined by \ $d=\mbox{min}\{ d(x,y):  x,y\in C, \ x\neq y\}$.
The \emph{weight} of a codeword $x$ is \ $w(x)=d(x,0)=|\{i  :  x_i\neq 0 \}|$.
For a linear code, $d=\mbox{min}\{ w(x)  :  x \in C, x\neq0 \}.$

A  $q$-ary linear code of length $n$, dimension $k$, and distance $d$ is called a $[n,k,d]_q$ code.  Either of the parameters $q$ and $d$ may be dropped when they are unspecified.
A linear $[n,k,d]$ code can detect at most $d-1$ errors in one codeword and correct at most \ $t = \left\lfloor  \frac{d-1}{2} \right\rfloor$ errors.
There are several upper bounds on $d$ for given parameters $n$ and $k$, we refer to \cite{Huffman} for numerous examples.
An $[n, k]$ linear code $C$ is \emph{optimal} if the minimum weight of $C$ achieves the theoretical upper bound on the minimum weight of $[n, k]$ linear codes.
An $[n, k]$ linear code $C$ is said to be a \emph{best known} linear $[n, k]$ code if $C$ has the highest minimum weight among all known $[n, k]$ linear codes.
A catalogue of best known codes is maintained in \cite{codetables}. We use this database to compare the minimum weight of all codes constructed in this paper.
Some stronger bounds are known for certain LCD codes (see for example some upper bounds on the minimum weight of binary LCD codes of small length computed via linear programming methods in \cite{ADDSS}), but over finite fields of order $q > 3$, due to a result of Carlet et al \cite[Corollary 5.3]{carlet-et-al}, the existence of a $[n,k,d]_{q}$ code implies the existence of a $[n,k,d]_{q}$ LCD code.

Let $w_{i}$ denote the number of codewords of weight $i$ in a linear code $C$ of length $n$.  Then the \emph{weight distribution} of $C$ is the list $[w_{i} : 0 \leq i \leq n]$.  A code is
\emph{formally self-dual} if $C$ and $C^{\perp}$ have identical weight distributions.  Formally self-dual even $[2n,n]$ codes may occasionally have a larger minimum weight than any $[2n,n]$
self-dual code \cite{KennedyPless}, and they may also be used to construct certain designs as a result of the Assmus-Mattson theorem.

Over fields of order $q^2$, the \emph{Hermitian dual} code $C^{\perp_H} = \{ v \in \F_{q^2}^n| \langle v,\overline{c}\rangle=0 {\rm \ for\ all\ } c \in C \}$, where
$\overline{c} = (c_{1}^q,\ldots,c_{n}^{q})$.  An alternative definition was introduced in \cite{Tonchev}.  There the \emph{Hermitian dual} code $C^{H}$ is the orthogonal complement of $C$ under the inverse dot product, given by $\langle u,v \rangle_{I} = u\cdot v^{*}$.  By \cite[Lemma 3.1]{WMat-paper}, this dot product is an inner product only when $q \in \{2,3,4\}$.
We adhere to the term Hermitian dual with respect to the inverse dot product when $q \in \{2,3,4\}$ as it coincides with the existing concept of Hermitian dual over $\F_{4}$,
and the use of the term \emph{Hermitian inner product} over $\F_{4}$ in \cite{Huffman}.

A linear code $C$ over $\F_{q^{2}}$ is usually called a Hermitian LCD code if $C\cap C^{\perp_H}=\{0\}$.  In the final section of this paper we adhere to the definition of the Hermitian dual $C^{H}$.  Note that $C^{H} = C^{\perp_H}$ when $q = 4$, and $C^{H} = C^{\perp}$ when $q \in \{2,3\}$.

\bigskip

In \cite{massey} Massey showed that the nearest-codeword (or maximum-likelihood) decoding problem for an LCD code may be simpler than that for a general linear code.
The decoding method proposed by Massey is based on the following statement.

\begin{thm}[Proposition 4, \cite{massey}]\label{massey}
Let $C$ be a LCD code of length $n$ over the field $\F_q$ and let $\varphi$ be a map $\varphi: C^\perp \rightarrow C$ such that $u \in C^\perp$ maps to one of the closest codewords $v$ to it in $C$.
Further, let $\Pi_C$ and $\Pi_{C^\perp}$ be the orthogonal projectors from $\F_q^n$ onto $C$ and $C^\perp$, respectively.
Then the map $\tilde{\varphi}: \F_q^n \rightarrow C$ such that
$$\tilde{\varphi}(w)= \Pi_C(w) + \varphi(\Pi_{C^\perp}(w))$$
maps each $w \in \F_q^n$ to one of it closest neighbours in $C$.
\end{thm}

\bigskip

For some of the codes constructed in this paper we define the map $\varphi$ partially and deduce a decoding algorithm, similarly to the approach employed in \cite{LCD-Peisert}.

\section{Euclidean LCD codes from weighing matrices}\label{lcd_mat}

The following two lemmas characterize Euclidean and Hermitian LCD codes, as shown in \cite{Carlet3,massey}.

\begin{lemma}[Proposition 1, \cite{massey}]\label{LCD-herm}
Let $G$ be a generator matrix for a code over a field. Then $\mathrm{det}(GG^{\top})\neq 0$ if and only if $G$ generates an LCD code.
\end{lemma}

Lemma \ref{LCD-herm} readily generalises to Hermitian codes over $\F_{q}$, $q \in \{2,3,4\}$.

\begin{lemma}[Proposition 2, \cite{Carlet3}]\label{LCD-herm2}
Let $G$ be a generator matrix for a Hermitian code over $\F_{q}$. Then $\mathrm{det}(GG^{*})\neq 0$ if and only if $G$ generates a Hermitian LCD code.
\end{lemma}

Let $O_{n}$ and $j_{n}$ denote the all zeros and all ones column vectors of length $n$.  The following simple result will be crucial so we state and prove it explicitly here.
\begin{prop}\label{det}
Let $A_1=aJ_n+xI_n$, $a,x\in \mathbb{F}$ for some field $\mathbb{F}$. The matrix $A_1$ is similar to the $(n\times n)$ matrix
$$A_2=\left[\begin{array}{c|c}
x+na & O_{n-1}^{\top} \\ \hline
aj_{n-1} & xI_{n-1} \end{array}\right].$$
Hence $\mathrm{det}(A_{1}) = (x+na)x^{n-1}$.
\end{prop}

{\bf Proof}.
Observe that $LA_{1}L^{-1} = A_{2}$ where
\[
L=\left[\begin{array}{c|c}
1 & j_{n-1}^{\top} \\ \hline
O_{n-1} & I_{n-1} \end{array}\right] ~~ \text{and} ~~ L^{-1}=\left[\begin{array}{c|c}
1 & -j^{\top} \\ \hline
O_{n-1} & I_{n-1} \end{array}\right].
\]
The result follows. {$\Box$}

\begin{thm}\label{des}
Let $W$ be a weighing matrix $\We(n,m)$, $B$ be the $n\times b$ point-by-block incidence matrix of a $(r,\lambda)$-design and let $G=\left[ \begin{array}{c|c}
W & B
\end{array} \right]$ be a matrix over $\F_{q}$.
If $(r+(n-1)\lambda+m)\neq 0$ and $(r-\lambda+m) \neq 0$ over $\F_{q}$, then $G$ generates a LCD code $C$ of length $n+b$ over $\F_{q}$.
\end{thm}

{\bf Proof}.
By construction $GG^{\top}=\lambda J_n+ (r-\lambda+m)I_n$ and using Proposition \ref{det} we have that $\mathrm{det}(GG^{\top}) = 0$ if and only if $(r+(n-1)\lambda+m) (r-\lambda+m)^{(n-1)}=0$ over
$\F_{q}$.  Thus by Lemma \ref{LCD-herm} $C$ is a LCD code if and only if $(r+(n-1)\lambda+m) (r-\lambda+m) \neq 0$ over $\F_{q}$. {$\Box$}

\bigskip

Choosing the appropriate $(r,\lambda)$-design allows us to control the dimension of the LCD code constructed.  As a special case, we may choose the design so that the incidence matrix is $I_{n}$.
This ensures that the code is a $[2n,n]$ code.  LCD codes can also be formally self-dual.  Certain $[2n,n]$ codes constructed in this paper will retain this property due to the following.

\begin{lemma}\label{weight dist}
Let $M$ be an $n \times n$ matrix over $\F_{q}$ and suppose $G=\left[ \begin{array}{c|c}
M & I_n
\end{array} \right]$ generates a $[2n,n]_{q}$ code $C$.  Further suppose that $\overline{G}=\left[ \begin{array}{c|c}
M & \alpha I_n
\end{array} \right]$ generates the dual code $C^{\perp}$ for some $\alpha \in \F_{q}$.  Then $C$ is formally self-dual.
\end{lemma}

{\bf Proof}.
Let $r_{1},\ldots,r_{n}$ and $\overline{r_{1}},\ldots,\overline{r_{n}}$ be the rows of $G$ and $\overline{G}$ respectively.  Then for any codeword $x = \textstyle{\sum}_{i=1}^{n}\lambda_{i}r_{i}$,
$\lambda_{i}\in \F_{q}$, the corresponding codeword $\overline{x} = \textstyle{\sum}_{i=1}^{n}\lambda_{i}\overline{r_{i}}$ is of the same weight.
Since the map $f : C \rightarrow C^{\perp}$ where $f(x) = \overline{x}$ for all $x \in C$ is bijective, the result follows.
{$\Box$}

\bigskip

Our goal is to construct LCD codes over various finite fields and for many parameters.  It is pertinent therefore to determine what families of matrices have sufficient properties so that as generator matrices, they satisfy the requirements of the lemmas above.  Weighing matrices and incidence matrices of $(r,\lambda)$-designs will most often provide the tools required.
Throughout this section $q=p^r$ is a prime power.  Working over $\F_{q}$, we will interpret integers as their value modulo $p$.  We consider Euclidean LCD codes throughout, until the final section.
As a result of Lemma \ref{weight dist}, certain $[2n,n]$ codes constructed will be formally self-dual. In all of these cases, the codes constructed are LCD provided a given condition
depending on the characteristic of $\F_{q}$ is met.  In the event that this condition is not met, the codes constructed are in fact self-dual.

\begin{cor}\label{detW}
Let $W$ be a weighing matrix $\We(n,m)$, and let $G=\left[ \begin{array}{c|c}
W & I_n
\end{array} \right]$ be a matrix over $\F_{q}$.
Then $GG^{\top}$ is nonsingular if and only if $m+1 \neq 0$ over $\F_{q}$.
\end{cor}

{\bf Proof}.
By construction $GG^{\top}=(m+1)I_{n}$ and so $\mathrm{det}(GG^{\top}) = 0$ if and only if $m+1 = 0$ over $\F_{q}$.
{$\Box$}

\begin{rem}
In the event that $m + 1 = 0$ over $\F_{q}$, the matrix $G$ above generates a self-dual $[2n,n]_{q}$ code.
\end{rem}

Using the characterization of LCD codes from \cite[Proposition 1]{massey} we have the following.

\begin{thm}\label{main}
Let $W$ be a weighing matrix $\We(n,m)$ and let $G=\left[ \begin{array}{c|c} W & I_n \end{array} \right]$ over $\F_{q}$. Then $G$ generates a $[2n,n]_{q}$ LCD code $C$ if $m + 1 \neq 0$ over $\F_{q}$.
\end{thm}

\begin{cor}\label{dual}
The matrix $\overline{G}=\left[ \begin{array}{c|c}
W & \alpha I_n
\end{array} \right]$ generates the dual code $C^{\perp}$ if $\alpha \in \F_{q}$ is such that $\alpha+m = 0$ over $\F_{q}$.
\end{cor}

{\bf Proof}.
Every row of $\overline{G}$ is orthogonal to every row of $G$, so the $n$-dimensional code generated by $\overline{G}$ is contained in $C^{\perp}$.
{$\Box$}

\bigskip

We denote the $ith$ row of the matrix $G$ by $r_i$, and the $ith$ row of $\overline{G}$ by $\overline{r}_i$ and prove the following lemma.

\begin{lemma}\label{lemma-decoding-1}
Let $W$ be a weighing matrix $\We(n,m)$ and $m + 1 \neq 0$ over the field $\F_{q}$. Further, let $C$ be the $[2n,n]_{q}$ LCD code generated by the matrix
$G=\left[ \begin{array}{c|c} W & I_n \end{array} \right]$. Suppose $C$ has minimum distance $d$ and $t=\left \lfloor \frac{d-1}{2} \right \rfloor$.
Then the following hold.
\begin{enumerate}
 \item Let $J$ be a subset of the set of indices of the rows of $W$. If $|J| \le t$, then the codeword in $C$ closest to $\sum_{i \in J} \lambda_i \overline{r}_i$,
       where $\lambda_i \in \F_{q}$ for $i \in J$, is $\sum_{i \in J} \lambda_i r_i$, at a distance $|J|$ from the vector $\sum_{i \in J} \lambda_i \overline{r}_i$.
 \item For $|J| \le t$ the map $\varphi$ of Theorem \ref{massey} can be uniquely defined by $\varphi(\sum_{i \in J} \lambda_i \overline{r}_i)=\sum_{i \in J} \lambda_i r_i$.
 \item If $w=\sum_{i \in J} \lambda_i \overline{r}_i$ and $w=\sum_{i \in K} \mu_i \overline{r}_i$, then $J=K$ and $\lambda_i=\mu_i$ for all $i \in J$.
\end{enumerate}
\end{lemma}

{\bf Proof}.
The proof of parts (1) and (2) follows from the definition of the matrices $G$ and $\overline{G}$ and the fact that $C$ has minimum distance $d$; and that of part (3) follows from the structure of the matrix $\overline{G}$.
{$\Box$}

\begin{lemma}\label{lemma-decoding-2}
Let $W$ be a weighing matrix $\We(n,m)$ and $m + 1 \neq 0$ over the field $\F_{q}$ and $C$ be the $[2n,n]_{q}$ LCD code generated by the matrix $G=\left[ \begin{array}{c|c} W & I_n \end{array} \right]$.
Suppose $C$ has minimum distance $d$ and $t=\left \lfloor \frac{d-1}{2} \right \rfloor$.
If the transmitted codeword from $C$ has no more than $t$ errors, it can be correctly decoded.
\end{lemma}

{\bf Proof}.
Suppose a codeword $c$ is sent and $w=c+w'$ is received, where $w'$ has no more than $t$ non-zero coordinates. Then $w'=\sum_{i \in K_1} \lambda'_i r_i + \sum_{i \in K_2} \mu'_i \overline{r}_i$,
where $|K_1|+|K_2| \le t$. Hence, $w=\sum_{i \in J_1} \lambda_i r_i + \sum_{i \in J_2} \mu_i \overline{r}_i= c+ \sum_{i \in K_1} \lambda'_i r_i + \sum_{i \in K_2} \mu'_i \overline{r}_i$.
It follows that
$\sum_{i \in J_1} \lambda_i r_i=c+ \sum_{i \in K_1} \lambda'_i r_i$ and $\sum_{i \in J_2} \mu_i \overline{r}_i=\sum_{i \in K_2} \mu'_i \overline{r}_i$. By Lemma \ref{lemma-decoding-1} (3),
$J_2=K_2$ and $\mu_i=\mu'_i$ for $i \in J_2$. Hence, $|J_2| \le t$ and $\tilde{\varphi}(w)=\sum_{i \in J_1} \lambda_i r_i + \varphi(\sum_{i \in J_2} \mu_i \overline{r}_i) =
\sum_{i \in J_1} \lambda_i r_i + \sum_{i \in J_2} \mu_i {r}_i$. Since $c=\sum_{i \in J_1} \lambda_i r_i - \sum_{i \in K_1} \lambda'_i r_i$, it follows that
$\tilde{\varphi}(w) - c=\sum_{i \in J_1} \lambda_i r_i + \sum_{i \in J_2} \mu_i {r}_i - \sum_{i \in J_1} \lambda_i r_i + \sum_{i \in K_1} \lambda'_i r_i=
\sum_{i \in J_2} \mu_i {r}_i + \sum_{i \in K_1} \lambda'_i r_i=\sum_{i \in K_2} \mu'_i {r}_i + \sum_{i \in K_1} \lambda'_i r_i$. Hence, $\tilde{\varphi}(w) - c$ is a codeword of $C$ that has
weight $|K_1|+|K_2| \le t$, so $\tilde{\varphi}(w) = c$.
{$\Box$}

\bigskip

Lemmas \ref{lemma-decoding-1} and \ref{lemma-decoding-2} show that the decoding method proposed by Massey when implemented for codes described in Corollary \ref{dual}, taking the map $\varphi$ from
Lemma \ref{lemma-decoding-1} (2), defined using the partial definition for linear combinations of at most $t=\left \lfloor \frac{d-1}{2} \right \rfloor$ rows of $\overline{G}$,
can correct up to $t$ errors.

\subsection{LCD codes from skew-weighing matrices}\label{skewMat}

A weighing matrix $\We(n,m)$ for which $W^{\top}=-W$ is called a skew-weighing matrix. It follows that a skew-weighing matrix has only zeros on its main diagonal. For further reading about skew-weighing matrices we refer the reader to \cite{seberry}.

\begin{thm}\label{code-skewMat}
Let $W$ be a skew-weighing matrix $\We(n,m)$, $B$ be the $n\times b$ point-by-block incidence matrix of a $(r,\lambda)$-design, and let $\alpha \in \F_{q}$. Then $G=\left[ \begin{array}{c|c}
W+\alpha I_n
&
B\\
\end{array} \right]$ generates a LCD code $C$ of length $n+b$ where $(m+\alpha^2+r+(n-1)\lambda)(m+\alpha^2+r-\lambda)^{n-1} \neq 0$ over the field $\F_{q}$.
\end{thm}

{\bf Proof}.
By construction $GG^{\top}=(m+\alpha^2+r - \lambda)I_n + \lambda J_{n}$ over $\F_{q}$.  Proposition \ref{det} and Lemma \ref{LCD-herm} complete the proof.
{$\Box$}

\begin{cor}\label{dual-skew}
Let $G$ be as above and suppose $B = I_{n}$.  Then $G$ generates a $[2n,n]_{q}$ LCD code $C$ where $m + \alpha^2 + 1 \neq 0$ over $\F_{q}$.  Further, the matrix $\overline{G}=\left[ \begin{array}{c|c}
W+\alpha I_n & \beta I_n\\
\end{array} \right]$, where $\beta + \alpha^2 + m = 0$ over $\F_{q}$, generates the dual code $C^{\perp}$.
\end{cor}

{\bf Proof}. The proof follows using arguments that are similar to those used in the proof of Corollary \ref{dual}. So we omit it. 
{$\Box$}

\bigskip

A Hadamard matrix $H$ is called skew-Hadamard if $H+H^{T}=2I$. The best known family of skew-Hadamard matrices are the Paley type I Hadamard matrices \cite{Paley} (see also Section \ref{ex-PaleyI} below).
In \cite{ha-mu} Harada and Munemasa used skew-Hadamard matrices of order 20 to classify self-dual $[20,10,9]_7$ codes. Theorem \ref{code-skewMat} can be applied to skew type Hadamard matrices, giving rise to the following results.

\begin{cor}\label{code-skewHadamard}
Let $H$ be a skew type Hadamard matrix. Then $G=\left[ \begin{array}{c|c}
H+\alpha I_n & B
\end{array} \right]$ generates a LCD code $C$ where $\alpha \in \F_{q}$ is such that $(n+\alpha^2 + 2\alpha + r + (n-1)\lambda)(n+\alpha^2 + 2\alpha + r -\lambda)^{n-1} \neq 0$ over $\F_{q}$.
\end{cor}

\begin{cor}\label{code-skewHadamard-I}
Let $G$ be as in Corollary~\ref{code-skewHadamard} and suppose $B = I_n$. Then $G$ generates a $[2n,n]_{q}$ LCD code $C$ when $n+(\alpha+1)^{2} \neq 0$ over $\F_{q}$.
Moreover, $\overline{G}=\left[ \begin{array}{c|c}
H+\alpha I_n & \beta I_n\\
\end{array} \right]$ where $\beta + (\alpha+1)^2 + n - 1 = 0$ over $\F_{q}$ generates the dual code $C^{\perp}$.
\end{cor}


\subsection{Examples from conference matrices and Paley type I Hadamard matrix} \label{ex-PaleyI}

In this section we illustrate with examples the constructions outlined in the previous sections. In the tables presenting codes over finite fields, $*$ denotes that the code is best known.

Recall that $O_{n}$ and $j_{n}$ denote the all zeros and all ones column vectors of length $n$.  Let $\pi$ be a prime power and let $\chi : F_{\pi} \rightarrow \{0,1,-1\}$ be the quadratic character where
$\chi(0) = 0$, $\chi(x) = 1$ if $x$ is a quadratic residue, and $\chi(x) = -1$ otherwise.
Let $A = [\chi(y-x)]_{x,y \in F_{\pi}}$. When $\pi \equiv 3 \;(\bmod\; 4)$ the Paley type I Hadamard matrix $\mathcal{P}_{1}(\pi)$ of order $\pi+1$ is
$$\mathcal{P}_{1}(\pi) = \left[ \begin{array}{cc} 1 & -j_{\pi}^{\top} \\ j_{\pi} & -A+I_{\pi} \end{array} \right].$$

Note that $\mathcal{P}_{1}(\pi)$ is skew, with all 1s on the main diagonal.

A conference matrix of order $t$ is a weighing matrix $W(t,t-1)$ with zeros on the main diagonal. If $\pi \equiv 1 \mod 4$ is a prime power, then the a similar construction yields a symmetric conference matrix of order $t=\pi+1$. For more information on conference matrices we refer the reader to \cite{crc-conf-mat}.

We first construct an $(r,\lambda)$-design to use for our examples.  Let $D_{1}$ be the point-by-block incidence matrix of a $\mathrm{BIBD}(15,35,14,6,5)$, and let $D_{2}$ be the point-by-block incidence matrix of a $\mathrm{PBD}(15,\{6,9\},6)$ with 35 blocks and 15 points, each of which being incident with 15 blocks.  Both $D_{1}$ and $D_{2}$ are obtained using a method of \cite{siamese}.  Then
$$B=\left[ \begin{array}{c|c}
j_{15} & D_{1}\\ \hline
O_{15} & D_{2}
\end{array} \right]$$ is the incidence matrix of an $(15,6)$-design with $30$ points and $36$ blocks.  Finally denote by $B_{t}$ the $t \times 36$ submatrix of $B$ comprising the first $t$ rows. Using Theorem \ref{des} we construct the matrices
$G_{t}=\left[ \begin{array}{c|c}
W_{t} & B_{t}
\end{array} \right]$ where $W_{t}$ is either a Paley type I Hadamard matrix or conference matrix of order $t$ as constructed above.  In Table \ref{table-des} we give examples of LCD codes obtained from $G_{t}$ according to whether $t \equiv 0,2 \;(\bmod\; 4)$.

\begin{table}[H]
\begin{center} \begin{footnotesize}
\begin{tabular}{|c|c|c||c|c|c|}
 \hline
$t \equiv 2 \;(\bmod\; 4)$  &  $C$ & Dual($C$) &$t \equiv 0 \;(\bmod\; 4)$  &  $C$& Dual($C$)  \\
\hline \hline
6 &  $[42,6,20]_3$ & $[42,36,2]_3$   &8& $[44,8,20]_3$& $[44,36,3]_3$   \\
10 &  $[46,10,22]_5$ & $[46,36,3]_5$  &  12& $[48,12,8]_2$ & $[48,36,2]_2$  \\
14 &  $[50,14,14]_3$ & $[50,36,5]_3$  &  12& $[48,12,18]_5$ & $[48,36,4]_5$  \\
14 &  $[50,14,18]_5$  & $[50,36,6]_5$    &  20& $[56,20,6]_2$ & $[56,36,2]_2$  \\
18 &  $[54,18,9]_3$  & $[54,36,6]_3$   &  20& $[56,20,8]_3$ & $[56,36,6]_3$  \\
18 &  $[54,18,18]_5$ & $[54,36,6]_5$     &  20& $[56,20,16]_5$ & $[56,36,7]_5$  \\
26 &  $[62,26,6]_3$ & $[62,36,6]_3$   &  24& $[60,24,4]_2$  & $[60,36,2]_2$ \\
30 &  $[66,30,5]_3$  & $[66,36,6]_3$  &&   &  \\
\hline \hline
\end{tabular} \end{footnotesize}
\caption{\footnotesize LCD codes constructed from weighing matrices and a $(15,6)$-design} \label{table-des}
\end{center}
\end{table}

In Table \ref{table-skewWeighing} we list LCD codes constructed by the method described in Section \ref{skewMat}, Corollary \ref{code-skewHadamard}, using the Paley type I skew-Hadamard matrices of orders up to 60 ($n=\pi+1=60$), and letting $B = I_{n}$.  These are $[2n,n]_{q}$ LCD codes, where $\alpha \in \F_{q}$ is such that $n+(\alpha+1)^2 \neq 0$. Recall that these are also formally self-dual codes.

\begin{table}[H]
\begin{center} \begin{footnotesize}
\begin{tabular}{|c | c  |c||c | c | c |}
 \hline
$n$  & $\alpha$& $C$ &  $n$  & $\alpha$& $C$ \\
\hline \hline
4 & 0 & $[8,4,2]_2$  &20 & 1& $[40,20,13]_5$*  \\
4 & 2 & $[8,4,3]_3$  & 24 & 0& $[48,24,2]_2$  \\
4 & 0 & $[8,4,4]_3$*  & 24 & 0& $[48,24,9]_3$  \\
4 & 1 & $[8,4,4]_5$*  & 24 & 1& $[48,24,15]_5$*  \\
8 & 0 & $[16,8,2]_2$    & 28 & 0& $[56,28,2]_2$  \\
8 & 2 &$[16,8,6]_3$*    & 28 & 2& $[56,28,6]_3$  \\
8 & 0& $[16,8,6]_5$    & 28 & 0& $[56,28,12]_3$  \\
8 & 1& $[16,8,7]_5$*   & 28 & 1& $[56,28,12]_5$  \\
12 & 0& $[24,12,2]_2$  & 28 & 0& $[56,28,15]_5$  \\
12 & 0& $[24,12,6]_3$   & 32 &  0& $[64,32,2]_2$  \\
12 & 1& $[24,12,6]_5$  & 32 & 2&  $[64,32,14]_3$  \\
12 &  0&$[24,12,8]_5$   & 32 & 2&  $[64,32,10]_5$  \\
12 &  4&$[24,12,9]_5$*  & 32 &  0& $[64,32,18]_5$*  \\
20 & 0& $[40,20,2]_2$  & 48 & 0& $[96,48,2]_2$  \\
20 & 2&$[40,20,10]_3$  & 48 & 0&$[96,48,15]_3$  \\
20 & 0&$[40,20,8]_5$  & 48 & 0& $[96,48]_5$ \\
\hline \hline
\end{tabular} \end{footnotesize}
\caption{\footnotesize $[2n,n]_{q}$ formally self-dual LCD codes constructed from skew-weighing matrices} \label{table-skewWeighing}
\end{center}
\end{table}

\begin{rem}
Some of the best known codes listed in Table \ref{table-skewWeighing} are optimal.
The optimal codes are $[8,4,4]_3$, $[8,4,4]_5$, $[16,8,6]_3$ and $[16,8,7]_5$ (see \cite{codetables}).
\end{rem}

\section{Euclidean LCD codes from orbit matrices of weighing matrices} \label{lcd_orb_mat}

Let $M$ be an $n \times n$ matrix with entries in some ring or field with additive identity $0$ and multiplicative identity $1$.  A \emph{permutation automorphism} of $M$ is a pair of permutation matrices $(P,Q)$ such that $PMQ^{\top} = M$.  The set of all such pairs forms a group under the composition $(P,Q)(R,S) = (PR,QS)$, denoted $\mathrm{PAut}(M)$.  Any subgroup $G \leq \mathrm{PAut}(M)$ acts as a permutation group on both the rows and columns of $M$.

\begin{lemma}\label{InvTranAut}
If $(P,Q)$ is a permutation automorphism of a nonsingular matrix $M$, then $(P,Q)$ is a permutation automorphism of $(M^{-1})^{\top}$.
\end{lemma}
{\bf Proof}.
$PMQ^{\top} = M \Rightarrow QM^{-1}P^{\top} = M^{-1} \Rightarrow P(M^{-1})^{\top}Q^{\top} = (M^{-1})^{\top}$.
{$\Box$}

\bigskip

Let $G$ be a permutation automorphism group of a nonsingular matrix $M=[m_{ij}]$, acting in $t$ orbits on the set of rows and the set of columns of $M$.
Let us denote the $G$-orbits on rows and columns of $M$ by ${\cal R}_{1},\ldots ,{\cal R}_{t}$ and ${\cal C}_{1},\ldots ,{\cal C}_{t}$, respectively,
and put $|{\cal R}_{i}|=\Omega _{i} $ and $ |{\cal C}_{i}|=\omega_{i}, \ i=1, \ldots , t$.  By Lemma \ref{InvTranAut}, the row and column orbits of $(M^{-1})^{\top}$ are identical, assuming the rows and columns are labelled identically.

Let $M_{ij}$ be the submatrix of $M$ consisting of the rows belonging to the row orbit ${\cal R}_{i}$ and the column belonging to ${\cal C}_{j}$.
We denote by $\Gamma_{ij}$ and $\gamma_{ij}$ the sum of a row and column of $M_{ij}$, respectively.
The numbers $\Gamma_{ij}$ and $\gamma_{ij}$ are well-defined, i.e. they do not depend on the choice of the row and the column,
because the sums of entries of any two rows (or columns) of $M_{ij}$ are equal. The $t \times t$ matrix $R=[\Gamma_{ij}]$
is called a \emph{row orbit matrix} of $M$ with respect to $G$.
The $t \times t$ matrix $C=[\gamma_{ij}]$ is called a \emph{column orbit matrix} of $M$ with respect to $G$.

\begin{lemma} \label{lemma-product}
Let $G$ be a permutation automorphism group of a nonsingular matrix $M=[m_{ij}]$ of order $n$,
and let ${\cal R}_{1},\ldots ,{\cal R}_{t}$ and ${\cal C}_{1},\ldots ,{\cal C}_{t}$
be the $G$-orbits on the rows and columns of the matrices $M$ and $N = (M^{-1})^{\top}=[n_{ij}]$.
Further, let $\Gamma_{ij}$ and $\gamma_{ij}$ be defined as above, for the row orbit matrix of $M$ and the column orbit matrix of $N$ respectively.
Then
$$\textstyle{\sum}_{j=1}^t \Gamma_{ij} \gamma_{sj}= \delta_{is},$$
where $\delta_{is}$ is the Kronecker delta.
\end{lemma}
{\bf Proof}.
Let $x$ be a row from the row orbit ${\cal R}_{i}$, and $y$ be a column from the column orbit ${\cal C}_{j}$. Then
\begin{align*}
\textstyle{\sum}_{j=1}^t \Gamma_{ij} \gamma_{sj}&=\textstyle{\sum}_{j=1}^t\big(\textstyle{\sum}_{z \in {\cal C}_{j}}m_{xz}\big)\big(\textstyle{\sum}_{w \in {\cal R}_{s}}n_{wy}\big)=
\textstyle{\sum}_{j=1}^t \textstyle{\sum}_{z \in {\cal C}_{j}} \textstyle{\sum}_{w \in {\cal R}_{s}} m_{xz} n_{wy}\\
&=\textstyle{\sum}_{j=1}^t \textstyle{\sum}_{z \in {\cal C}_{j}} \textstyle{\sum}_{w \in {\cal R}_{s}} m_{xz} n_{wz}
=\textstyle{\sum}_{j=1}^t \textstyle{\sum}_{w \in {\cal R}_{s}} \textstyle{\sum}_{z \in {\cal C}_{j}} m_{xz} n_{wz}\\
&=\textstyle{\sum}_{w \in {\cal R}_{s}} \textstyle{\sum}_{j=1}^t \textstyle{\sum}_{z \in {\cal C}_{j}} m_{xz} n_{wz}=
\textstyle{\sum}_{w \in {\cal R}_{s}} \textstyle{\sum}_{z=1}^n m_{xz} n_{wz}.
\end{align*}
If $i \neq s$, then
$$\textstyle{\sum}_{w \in {\cal R}_{s}} \textstyle{\sum}_{z=1}^n m_{xz} n_{wz}=\textstyle{\sum}_{w \in {\cal R}_{s}} 0=0.$$
If $i = s$, then
$$\textstyle{\sum}_{w \in {\cal R}_{s}} \textstyle{\sum}_{z=1}^n m_{xz} n_{wz}= (\Omega_s-1)0+1=1,$$
where $\Omega_s$ is the length of the orbit ${\cal R}_{s}$.
{$\Box$}

\bigskip

\begin{thm}\label{R nonsingular}
If $R$ is the row orbit matrix of a nonsingular matrix $M$ with respect to a permutation automorphism group $G$, then $R$ is nonsingular. Moreover, $R^{-1}$ is the transpose of the column orbit matrix of $(M^{-1})^{\top}$ with respect to $G$.
\end{thm}

Construction of self-orthogonal codes from orbit matrices of block designs was introduced
in \cite{Harada-Tonchev} and further developed in \cite{DC-BGR-SR-LS}. In \cite{DMpaper} the authors defined orbit matrices of Hadamard
matrices and showed that under certain conditions the orbit matrices yield self-orthogonal
codes. In this section we give a construction of LCD codes using orbit matrices of weighing matrices.

The following result, originally presented in \cite{DMpaper}, was applied to Hadamard matrices. Here it has been modified for weighing matrices.  The proof is similar.

\begin{thm} \label{thm-product}
Let $G$ be a permutation automorphism group of a $\We(n,m)$ weighing matrix $W$ acting with $t$ orbits, then
$$\textstyle{\sum}_{j=1}^t \frac{\Omega_s}{\omega_j}\Gamma_{ij} \Gamma_{sj}= \delta_{is} m,$$
where $\delta_{is}$ is the Kronecker delta.
\end{thm}

In other words, the orthogonality of distinct rows of $W$ implies the orthogonality of distinct rows of $R$.

\begin{cor}
If $G$ is a permutation automorphism group of a weighing matrix $W$, acting with all orbits of the same size, then the row and column orbit matrices of $W$ with respect to $G$ are equal.
\end{cor}

{\bf Proof}.
By Theorem \ref{R nonsingular}, $R$ is of rank $t$ with inverse $m^{-1}R^{\top}$ where $m^{-1}R$ is the column orbit matrix of $(W^{-1})^{\top} = m^{-1}W$.  Thus $R$ is also the column orbit matrix of $W$ with respect to $G$.
{$\Box$}

\begin{thm} \label{thm-orbW}
Let $W$ be a weighing matrix $\We(n,m)$ and $G$ be a permutation automorphism group of $W$ acting with $t$ orbits, each of the same length.  Further let $B$ be the $t\times b$ point-by-block incidence matrix of a $(r,\lambda)$-design. If $R$ is the row orbit matrix of $W$ with respect to $G$, then $A=\left[ \begin{array}{c|c}
R & B
\end{array} \right]$ generates a LCD code $C$ of length $t+b$ over $\F_{q}$ where $(r + (t-1)\lambda + m) \neq 0$ and $(r -\lambda + m) \neq 0$ over $\F_{q}$.
\end{thm}

{\bf Proof}. Follows by using arguments that are similar to those used in the proof of Theorem \ref{des}. {$\Box$}

\begin{cor}\label{cor-orbW}
Let $A$ be the matrix of Theorem \ref{thm-orbW} and suppose that $B = I_{t}$.  Then $A$ generates a formally self-dual $[2t,t]_{q}$ LCD code where $m+1 \neq 0$ over $\F_{q}$. Moreover, $\overline{A}=\left[ \begin{array}{c|c}
R & \alpha I_t
\end{array} \right]$ generates the dual code $C^{\perp}$ where $\alpha + m = 0$ over $\F_{q}$.
\end{cor}

{\bf Proof}.
By Theorem \ref{thm-product}, $AA^{\top} = (m+1)I_{t}$ so the first claim follows from Lemma \ref{LCD-herm}.  The orthogonality of the rows of $A$ and $\overline{A}$ proves the latter.
{$\Box$}

\subsection{Examples}

In this subsection we construct LCD codes from various weighing matrices and their orbit matrices, adhering to the theory outlined in this section.  In particular we use Bush-type Hadamard matrices $H_{36}$ and $H_{100}$ of order 36 \cite{janko36} and of order 100 \cite{janko-hadi-tonchev}, and various weighing matrices constructed from periodic ternary Golay pairs via the construction described in \cite[Section 6]{WMat-paper}.  These examples were chosen for their automorphism groups satisfying the necessary orbit length properties of Theorem \ref{thm-orbW}.  In particular, weighing matrices of order $2n$ constructed from periodic ternary Golay pairs of length $n$ have cyclic permutation automorphism groups of order $k$ acting on rows and columns in $\frac{2n}{k}$ orbits of length $k$, for any $k$ dividing $n$.  Some of the ternary Golay pairs used are taken from \cite{Craigen}.

In Table \ref{table-had36} we give the parameters of the LCD codes constructed from orbit matrices of the matrices $H_{36}$ and $H_{100}$ with respect to subgroups of the permutation automorphism groups of $H_{36}$ and $H_{100}$, and applying Corollary \ref{cor-orbW}.

\begin{table}[H]
\begin{center} \begin{footnotesize}
\begin{tabular}{| c | c | c | c | c | c |}
 \hline
$G \leq \mathrm{PAut}(H_{36})$& $C$ & $G \leq \mathrm{PAut}(H_{36})$ & $C$ & $G \leq \mathrm{PAut}(H_{100})$ & $C$ \\
\hline \hline
$I$ & $[72,36,2]_2$ & $Z_3$ & $[24,12,2]_2$ & $Z_5$ & $[40,20,2]_2$ \\
$I$ & $[72,36,6]_3$ & $Z_3$ & $[24,12,3]_3$ & $Z_5$ & $[40,20,4]_3$ \\
$I$ & $[72,36,12]_5$ & $Z_3$ & $[24,12,8]_5$ & $Z_5$ & $[40,20,2]_5$ \\
$I$ & $[72,36,12]_7$ & $Z_3$ & $[24,12,8]_7$ & $Z_5$ & $[40,20,4]_7$ \\
$I$ & $[72,36,6]_9$ & $Z_3$ & $[24,12,3]_9$ & $Z_5$ & $[40,20,4]_9$ \\
$I$ & $[72,36,12]_{11}$ & $Z_3$ & $[24,12,8]_{11}$ & $Z_5$ & $[40,20,4]_{11}$ \\
$I$ & $[72,36,12]_{25}$ & $Z_3$ & $[24,12,8]_{25}$ & $Z_5$ & $[40,20,2]_{25}$ \\
\hline \hline
\end{tabular} \end{footnotesize}
\caption{\footnotesize Formally self-dual LCD codes obtained from Bush-type Hadamard matrices} \label{table-had36}
\end{center}
\end{table}

The codes in Table \ref{table-W72} are constructed from weighing matrices $\We(42,26)$, $\We(50,29)$, $\We(56,29)$ and $\We(72,36)$ which were constructed from periodic ternary Golay pairs.
In Table \ref{table-W72} we give the parameters of the LCD codes constructed from orbit matrices of these matrices, with respect to subgroups of the permutation automorphism groups, again using Corollary \ref{cor-orbW}.

\begin{table}[H]
\begin{center} \begin{footnotesize}
\begin{tabular}{| c | c | c | c | c | c |}
 \hline
$G \leq \mathrm{PAut}(W)$ & $C$ & $G \leq \mathrm{PAut}(W)$ & $C$ & $G \leq \mathrm{PAut}(W)$ & $C$ \\
\hline \hline
$I$ & $[144,72,12]_5$ &     $Z_4$ & $[36,18,6]_5$ &     $Z_5$ & $[20,10,6]_{11}$ \\
$I$ & $[144,72,12]_7$ &        $Z_4$ & $[36,18,6]_7$ &     $Z_7$ & $[16,8,6]_7$ \\
$I$ & $[144,72,12]_{11}$ &     $Z_4$ & $[36,18,6]_{11}$ &  $Z_7$ & $[16,8,6]_{11}$ \\
$I$ & $[144,72,12]_{25}$ &     $Z_4$ & $[36,18,6]_{25}$ &  $Z_9$ & $[16,8,4]_5$ \\
$Z_2$ & $[72,36,6]_5$ &     $Z_3$ & $[28,14,4]_2$ &     $Z_9$ & $[16,8,4]_7$ \\
$Z_2$ & $[72,36,6]_7$ &     $Z_3$ & $[28,14,8]_5$ &     $Z_9$ & $[16,8,4]_{11}$ \\
$Z_2$ & $[72,36,6]_{11}$ &  $Z_3$ & $[28,14,8]_7$ &     $Z_9$ & $[16,8,4]_{25}$ \\
$Z_2$ & $[72,36,6]_{25}$ &  $Z_3$ & $[28,14,10]_{11}$ & $Z_7$ & $[12,6,2]_2$ \\
$Z_3$ & $[48,24,4]_5$ &     $Z_3$ & $[28,14,8]_{25}$ &  $Z_7$ & $[12,6,4]_5$ \\
$Z_3$ & $[48,24,4]_7$ &     $Z_4$ & $[28,14,9]_7$ &     $Z_7$ & $[12,6,6]_7$ \\
$Z_3$ & $[48,24,4]_{11}$ &  $Z_4$ & $[28,14,10]_{11}$ & $Z_7$ & $[12,6,6]_{11}$ \\
$Z_3$ & $[48,24,4]_{25}$ &  $Z_5$ & $[20,10,6]_7$ &     $Z_7$ & $[12,6,4]_{25}$ \\
\hline \hline
\end{tabular} \end{footnotesize}
\caption{\footnotesize Formally self-dual LCD codes from orbit matrices of weighing matrices} \label{table-W72}
\end{center}
\end{table}

\subsection{Euclidean LCD codes from orbit matrices of skew-weighing matrices}\label{skewOM}

In Section \ref{skewMat} we showed how LCD codes can be constructed using skew-weighing matrices. In this subsection we extend the construction of LCD codes to orbit matrices of skew-weighing matrices.  Let $W$ be a skew-weighing matrix $\We(n,m)$ with rows and columns labelled $r_{1},\ldots,r_{n}$ and $c_{1},\ldots,c_{n}$ respectively. The following results are obtained.

\begin{lemma}
Let $G$ be a permutation automorphism group of $W$ acting with $t$ orbits of the same length $w$, such that $r_{s} \in \mathcal{R}_{i}$ if and only if $c_{s} \in \mathcal{C}_{i}$.  Then the orbit matrix $R$ is skew symmetric.
\end{lemma}

{\bf Proof}.
As the sum of the entries in the submatrix $W_{ij}$ is invariant whether we count row sums or column sums, it follows that $\Gamma_{ij} = \gamma_{ij}$ for all $1 \leq i,j \leq t$. As $W$ is skew, $\Gamma_{ij} = -\gamma_{ji} = -\Gamma_{ji}$, and consequentially $\Gamma_{ii} = 0$.
{$\Box$}

\begin{thm} \label{code-skew}
Let $G$ be a permutation automorphism group of $W$ acting with $t$ orbits of the same length $w$, such that $r_{s} \in \mathcal{R}_{i}$ if and only if $c_{s} \in \mathcal{C}_{i}$.  Further let $B$ be the $t\times b$ point-by-block incidence matrix of a $(r,\lambda)$-design. If $R$ is the row orbit matrix of $W$ with respect to $G$, then $A=\left[ \begin{array}{c|c}
R+\alpha I_t & B
\end{array} \right]$ generates a LCD code $C$ of length $t+b$ over $\F_{q}$ where $(r + \alpha^2 + (t-1)\lambda + m) \neq 0$ and $(r + \alpha^2 - \lambda + m) \neq 0$ over $\F_{q}$.
\end{thm}

\begin{cor}\label{dual-skewOM}
Let $A$ be the matrix of Theorem \ref{code-skew} and suppose $B = I_{t}$.  Then $A$ generates a formally self-dual $[2t,t]_{q}$ LCD code $C$, where $\alpha \in \F_{q}$ satisfies $m+\alpha^2+1 \neq 0$.  Moreover, $\overline{A}=\left[ \begin{array}{c|c}
R+\alpha I_t & \beta I_t
\end{array} \right]$ generates the dual code $C^{\perp}$ where $\beta +\alpha^2+m = 0$ over $\F_{q}$.
\end{cor}

Any permutation automorphism of a skew-weighing matrix $\We(n,n-1)$ preserves the main diagonal, and thus if a group $G$ acts with all orbits of equal length, it satisfies the conditions of Theorem \ref{code-skew}.

\section{Hermitian LCD codes from weighing matrices}\label{LCD-orbitMat2}

We conclude by generalizing the methods discussed so far to construct Hermitian LCD codes from $\F_{q}$-weighing matrices and their orbit matrices, where $q \in \{2,3,4\}$.  To demonstrate we obtain Hermitian LCD codes over $\F_{4}$, from a $\We(6,6;\F_{4})$, a $\We(8,8;\F_{4})$ and a $\We(12,12;\F_{4})$. An $\F_{q}$-weighing matrix $\We(n,m;\F_{q})$, introduced in \cite{WMat-paper},
is an $n \times n$ matrix $W$ with $m$ non-zero entries from $\F_{q}$ in each row and column such that $WW^{*} = mI_{n}$. In Table \ref{table-lcd-F4} we list the obtained results.

\begin{table}[H]
\begin{center} \begin{footnotesize}
\begin{tabular}{|c | c | c | c |}
 \hline
$n$  & $C$ & $\mathrm{Dual}(C)$\\
\hline \hline
6  & $[12,6,4]_4$  &$[12,6,4]_4$ \\
8  & $[16,8,4]_4$  &$[16,8,4]_4$ \\
12 &  $[24,12,4]_4$  &$[24,12,4]_4$\\
\hline \hline
\end{tabular} \end{footnotesize}
\caption{\footnotesize Hermitian LCD codes over $\F_{4}$} \label{table-lcd-F4}
\end{center}
\end{table}

The orbit matrix of an $\F_{q}$-weighing matrix is defined similarly to the orbit matrix of a weighing matrix, but with all arithmetic carried out over $\F_{q}$; see \cite{WMat-paper} for further details.

\begin{thm} \label{thm-orbWH}
Let $W$ be a $\We(n,m;\F_{q})$ where $q \in \{2,3,4\}$ and let $G$ be a permutation automorphism group of $W$ acting with $t$ orbits, each of the same length.  Further let $B$ be the $t\times b$ point-by-block incidence matrix of a $(r,\lambda)$-design. If $R$ is the row orbit matrix of $W$ with respect to $G$, then $A=\left[ \begin{array}{c|c}
R & B
\end{array} \right]$ generates a Hermitian LCD code $C$ of length $t+b$ over $\F_{q}$ where $(r + (t-1)\lambda + m) \neq 0$ and $(r -\lambda + m) \neq 0$ over $\F_{q}$.
\end{thm}

\begin{cor}\label{cor-orbWH}
Let $A$ be the matrix of Theorem \ref{thm-orbWH} and suppose that $B = I_{t}$.  Then $A$ generates a $[2t,t]_{q}$ Hermitian LCD code where $m+1 \neq 0$ over $\F_{q}$. Moreover, $\overline{A}=\left[ \begin{array}{c|c}
R & \alpha I_t
\end{array} \right]$ generates the dual code $C^{\perp}$ where $\alpha + m = 0$ over $\F_{q}$.
\end{cor}

\begin{rem}
In Lemmas \ref{lemma-decoding-1} and \ref{lemma-decoding-2} we showed that the decoding method proposed by Massey can be implemented for codes described in Corollary \ref{dual}.
Codes obtained by using Corollaries \ref{code-skewHadamard-I}, \ref{cor-orbW}, \ref{dual-skewOM} and \ref{cor-orbWH} also allow decoding using the decoding method proposed by Massey, where
the map $\varphi$ can be defined as in Lemma \ref{lemma-decoding-1} (2).
\end{rem}


\noindent {\bf Acknowledgement} \\
D. Crnkovi\' c and A. \v Svob were supported by {\rm C}roatian Science Foundation under the project 6732.  R. Egan was supported by the Irish Research Council (Government of Ireland Postdoctoral Fellowship, GOIPD/2018/304). B. G. Rodrigues work is based on the research supported by the National Research Foundation of South Africa (Grant Numbers 95725 and 106071). B. G. Rodrigues acknowledges support from the Erasmus Mundus Plus academic exchange programme to visit the University of Rijeka in 2018.

\end{document}